\documentclass[amstex,11pt,a4paper]{article}

\usepackage[version=0.96]{pgf}
\usepackage{tikz}
\usetikzlibrary{arrows,shapes,snakes,automata,backgrounds,petri}

\usepackage[english,russian]{babel}
\usepackage{amsfonts,amssymb}
\usepackage{graphicx}
\usepackage{latexsym,exscale}
\usepackage{amsmath}
\usepackage{amsthm}
\usepackage{blindtext}

\usepackage{amsmath}
\usepackage{relsize}

\usepackage{xcolor}
\usepackage[hypertex,linkbordercolor={0 1 1}]{hyperref}

\usepackage[version=0.96]{pgf}
\usepackage{tikz}
\usetikzlibrary{arrows,shapes,snakes,automata,backgrounds,petri}
\usepackage[latin1]{inputenc}
\usetikzlibrary{intersections}

\definecolor{linkcolor}{HTML}{F70A0A}
\definecolor{urlcolor}{HTML}{F70A0A}

\newcommand{\bigstart}{\mathop{\Large \mathlarger{\mathlarger{*}}}}

\newtheorem{Lemma}{Lemma}
\newtheorem{Theorem}{Theorem}
\newtheorem*{myth}{Theorem}
\newtheorem{Corollary}{Corollary}
\newtheorem{Definition}{Definition}

\newtheorem{Example}{Example}

\newenvironment{Proof} 
{\par\noindent{\bf{Proof.}}} 
{\hfill$\scriptstyle\blacksquare \\$}

\newcounter{Remark}

\renewcommand{\ge}{\geqslant}

\title{On free decompositions of verbally closed subgroups in free products of finite groups}
\author{A. M. Mazhuga \\ email \href{mailto:mazhuga.andrew@yandex.ru}{mazhuga.andrew@yandex.ru}
   }
\date{}

\makeatletter
\let\newtitle\@title
\let\newauthor\@author
\let\newdate\@date
\makeatother

\begin{document}

\maketitle

\begin{otherlanguage}{english}

\begin{abstract}

For a subgroup of a free product of finite groups, we obtain necessary conditions (on its  Kurosh decomposition) to be verbally closed.
\end{abstract}

\end{otherlanguage}

\section{Introduction}

Algebraically closed objects play important part in modern algebra. In the present paper, we study verbally closed subgroups in free products of finite groups. First of all we establish some terminology.

For groups $H$ and $G$ we write $H\leqslant G$ if $H$ is a subgroup of $G$. Let $X=\left \{ x_i\ |\ i\in\mathbb{N} \right \}$ be a countable infinite set of variables, and $F(X)$ be the free group with basis $X$. An \emph{equation} with variables $x_1,\dots, x_n\in X$ and coefficients from $H$ is an arbitrary expression of the form $w(x_1,\dots,x_n,H)=1$, where $w(x_1,\dots,x_n,H)$ is a word in the alphabet $X\cup  X^{-1}\cup H$ (in other words, $w(x_1,\dots,x_n,H$) lies in the free product $F(X)\ast H$ of $F(X)$ and $H$). If the left-hand side of the equation does not depend on coefficients from $H$ we will omit $H$ in its expression: $w(x_1,\dots,x_n)=1$. We say that an equation $w(x_1,\dots,x_n,H)=1$ has a \emph{solution} in $G$ if there are some elements $g_i\in G$, $i=1,\dots,n$ such that $w(g_1,\dots,g_n,H)=1$ in $G$.

A subgroup $H$ is called \emph{algebraically closed} in $G$ if for every finite system of equations $S = \left \{ w_i\left ( x_1,\dots,x_n,H \right )=1\ |\ i=1,\dots, m \right \}$ with coefficients from $H$ the following holds: if $S$ has a solution in $G$, then it has a solution in $H$. The following notion of verbally closed subgroups was introduced by V.~ Roman'kov and A.~Myasnikov in~\cite{Romonkov}.

{\Definition\ }\label{VCS}\!\!\! A subgroup $H$ of $G$ is called \emph{verbally closed} if for any $w\in F(X)$ and $h\in H$ an equation $w(x_1,\dots,x_n)h^{-1} = 1$ has a solution in $H$ if it has a solution in $G$.

Verbally closed subgroups have a connection to research on verbal width of elements in groups. For $w=w(x_1,\dots,x_n)\in F(X)$ and a group $G$ by $w[G]$ we denote the set of all $w$-elements in $G$, i.e. $w[G] = \left \{ w(g_1,\dots,g_n)\ |\ g_1,\dots,g_n\in G \right \}$. The \emph{verbal subgroup} $w(G)$ is the subgroup of $G$ generated by $w[G]$. Therefore $H$ is a verbally closed subgroup in $G$ if and only if $w[H] = w[G]\cap H$ for every $w\in F(X)$. The \emph{$w$-width} $l_{w,G}(g)$ of an element $g\in w(G)$ is the minimal natural number $n$ such that $g$ is a product of $n$ $w$-elements or their inverses in $G$; the width of $w(G)$ is the supremum of widths of its elements. Usually, it is very hard to compute the $w$-length of a given element $g\in w(G)$ or the width of $w(G)$. For more details we refer to D.~Segal's book~\cite{Segal}.

Another important notion for us is a rectact.

{\Definition\ }\!\!\! A subgroup $H$ of a group $G$ is called a \emph{retract} of $G$ if there is a homomorphism (termed \emph{retraction}) $\varphi\colon G \to H$ which is identical on $H$. Equivalently, a retraction $\varphi\colon G \to G$ is a homomorphism such that $\varphi^2 = \varphi$ and a retract is the image of a retraction.

Groups algebraically closed in the class of all groups (a group $H$ is algebraically closed in a class of groups $\mathfrak{K}$ if and only if $H$ is algebraically closed in every extension $H\leqslant G$ with $G\in \mathfrak{K}$) were introduced by W.~R.~Scott in~\cite{Scott}. They have been thoroughly studied in 1970-80's, see for references article~\cite{Macintyre} and books~\cite{Higman},~\cite{Hodges}. On the other hand, we do not know much about algebraic closedness of subgroups in particular groups.

A several papers devoted to study of retracts of groups. It is known that the class of retracts of a non-abelian free group is much larger than its class of free factors (see, for instance,~\cite{Martino}). In~\cite{Bergman} J.~Bergman proved that intersection of finitely many of retracts of a free group is again a retract. At the same time, the question remains open whether it is true that if $F$ is a free group of finite rank, $R$ a retract of $F$, and $H$ a subgroup of $F$ of finite rank, then the intersection $H\cap R$ is a retract of $H$ (see Question~19 in~\cite{Bergman}). An element $g\in G$ is called a \emph{test element} (for recognizing automorphisms) of a group $G$ if, whenever $\varphi(g) = g$ for an endomorphism $\varphi$ of $G$, it follows that $\varphi$  is an automorphism. In~\cite{Turner} E.~C.~Turner obtained the following result and characterized test elements of a free group as elements that do not belong to any proper retract.

\begin{myth}
If $\varphi\colon F \to F$ is an endomorphism of the finitely generated free group $F$, then the subgroup

$$H_\varphi  = \bigcap_{n=0}^{\infty} \varphi^{n}\left ( F \right )$$

is a retract of $F$; $H_\varphi$ is a proper retract precisely when $\varphi$ fails to be an automorphism. If $\varphi$ is a monomorphism, then $H_\varphi$ is a free factor.
\end{myth}
The generalization of this theorem for hyperbolic groups was proved in~\cite{Neill}.

Very little is known in general about verbally closed subgroups of a given group $G$. For instance, the following questions are still open for most non-abelian groups:

\begin{enumerate}
    \item[\sffamily $\mathrm{1.}$] Is there an algebraic description of verbally closed subgroups of $G$?

    \item[\sffamily $\mathrm{2.}$] Is the class of verbally closed subgroups of $G$ closed under intersections?

    \item[\sffamily $\mathrm{2.}$] Does there exist the verbal closure $vcl(H)$ ($vcl(H)$ is least (relative to inclusion) verbally closed subgroup of $G$ containing $H$) of a given subgroup $H$ of $G$?
    \end{enumerate}

Now we shall briefly discuss the relationship among retracts, algebraically closed and verbally closed subgroups of groups.  In general, the property of being verbally closed is essentially weaker than being algebraically closed (see examples of verbally but not algebraically closed subgroups in 2-nilpotent torsion-free groups in~\cite{Baumslag}). It is easy to prove (see, for instance,~\cite{Romonkov}) that, if $H$ is a retract of $G$, then $H$ is algebraically (and verbally) closed in $G$. If $G$ is a finitely presented group, and $H$ is a finitely generated subgroup of $G$, then $H$ is an algebraically closed subgroup of $G$ if and only if $H$ is a retract of $G$ (see~\cite{Romonkov}). It was proved in~\cite{Romonkov} and~\cite{Romonkov2} that a subgroup $H$ of a free or free nilpotent group $G$ of finite rank is verbally closed if and only if $H$ is a retract of $G$. This implies that the properties of being verbally closed, being algebraically closed, and being a retract are equivalent for subgroups of free or free nilpotent groups of finite rank. Also in~\cite{Romonkov} the following result was obtained.

\begin{myth}
 If $F$ is a free group of finite rank, then the following holds:

  \begin{enumerate}
    \item[\sffamily $\mathrm{1.}$] There is an algorithm to decide if a given finitely generated subgroup of $F$ is verbally (algebraically) closed or not.

    \item[\sffamily $\mathrm{2.}$] There is an algorithm to find a basis of $vcl(H)$ for a given finitely generated subgroup $H$ of $F$.

  \end{enumerate}

 \end{myth}

Except for free, some nilpotent and finite groups the structure of verbally closed subgroups remains unknown. The main problem in study of verbal closeness is that it is exceedingly difficult to apply Definition~\ref{VCS} straightforwardly.

The following theorem was obtained by A.~Kurosh~\cite{Kurosh1}.

  \begin{myth}\label{Kur}
 Let $$G = \bigstart_{\alpha \in \Lambda }A_\alpha$$ be a free product of groups $A_\alpha$ and $H$ be a subgroup of $G$. Then $H$ is a free product $$H = F\ast\bigstart_{\beta}B_\beta,$$ where $F$ is a free group such that $F\cap A_\alpha^g = \{1\}$ for all $\alpha \in \Lambda $, $g\in G$ and each subgroup $B_\beta$ is conjugate in $G$ to a subgroup of $A_\alpha$ for some $\alpha$.
 \end{myth}

Such free decompositions of a subgroup $H$ of $G$ we call \emph{Kurosh free decompositions} (of a subgroup $H$).

The main goal of this article is to study structure of verbally closed subgroups of free products of finitely many finite groups. For a subgroup of this type of groups we obtain necessary conditions (on its  Kurosh decomposition) to be verbally closed (Theorem~\ref{ThM}). To prove the first part of Theorem~\ref{ThM} we use an elementary technique combined with some tricks (Lemma~\ref{LemF}). To prove its second part we use a bit of Bass-Serre theory (Lemma~\ref{BFL}). To prove that this conditions are not sufficient we use Lee's $C$-test words (the proof of Theorem~\ref{ThLee1}).

Now we fix some notations that will be used in this text.

Denote by $h^{gn}$, $n\in\mathbb{Z}$ an element $gh^ng^{-1}$ (similarly, $H^g=gHg^{-1}$), and by $[g,h]$ a commutator of the form $ghg^{-1}h^{-1}.$ A cyclic group of order $n$ generated by an element $a$ will be denoted as $\left \langle a \right \rangle_n$, i.e. $\left \langle a \right \rangle_n\simeq \mathbb{Z}/n\mathbb{Z}$ (similarly, $\left \langle a \right \rangle_\infty\simeq \mathbb{Z}$). $\#G$ is the cardinality of a set $G$. Sometimes we use an abridged notation of the type $\underline{x}=x_1,\dots,x_n$ or even $\underline{x_i'} = (x_{i,1}',\dots,x_{i,n}')$ rewriting, for example, an equation $w(x_1,\dots,x_n)h^{-1} = 1$ as $w(\underline{x})h^{-1} = 1$.

\section{The main result}

The following theorem is the main result of this article.

{\Theorem\label{ThM} Let $H$ be a subgroup of $G=G_1\ast \cdots\ast G_n$, $1<\#G_i<\infty$, $i=1,\dots,n$ with a Kurosh free decomposition:

\begin{equation}\label{Dec}
H = F\ast\bigstart_{i=1}^{n}\left(\bigstart_{j}^{}H_{i,j}^{g_{i,j}}\right),
\end{equation}
where $F$ is a free group and $H_{i,j}\ne \left \{ 1 \right \}$ is a subgroup of $G_i$ for all $j$.

If $H$ is a verbally closed subgroup in $G$, then the following holds:

\begin{enumerate}
    \item[\sffamily $\mathrm{1.}$] $F$ is the trivial group.

    \item[\sffamily $\mathrm{2.}$] For any two subgroups $H_{i,j_1}$ and $H_{i,j_2}$ from decomposition (\ref{Dec}) we have: if there are elements $f_i$, $g_i \in G_i$ and natural numbers $k_1$, $k_2$ such that $f_i^{k_1}\in H_{i,j_1}$ and $f_i^{k_2}\in H_{i,j_2}^{g_i}$, then either $f_i^{k_1}=1$ or  $f_i^{k_2}=1$.
    \end{enumerate}
}

The second assertion of the theorem means that subgroups $H_{i,j_1}$ and $H_{i,j_2}^{g_i}$ cannot simultaneously intersect the same cyclic subgroup nontrivially. For example, the following condition can be easily deduced from the second assertion of the theorem:

\begin{enumerate}
    \item[\sffamily $\mathrm{3.}$] For any two subgroups $H_{i,j_1}$ and $H_{i,j_2}$ from decomposition (\ref{Dec}) we have: for an element $g_i\in G_i$ the intersection $H_{i,j_1}\cap H_{i,j_2}^{g_i}$ is trivial.
    \end{enumerate}

The following examples and corollaries demonstrates how the second assertion of this theorem can be applied.

{\Example\ }\!\!\!\!\! Let $H = \left \langle a \right \rangle_2\ast\left \langle b^c \right \rangle_2$ be a subgroup of $G=G_1\ast G_2 = \left \langle a,b\ \!|\ \!a^2,b^2,(ab)^3 \right \rangle\ast\left \langle c\ \!|\ \!c^2 \right \rangle\simeq S_3\ast\mathbb{Z}_2$. Then $H$ is not a verbally closed subgroup in $G$ since for $g_1=ba\in G_1$ we have $\left \langle a \right \rangle_2\cap\left \langle b^{g_1} \right \rangle_2 = \left \langle a \right \rangle_2\cap\left \langle a \right \rangle_2\ne \{1\}$ and it is the contradiction to the second assertion of the theorem.

{\Example\ }\!\!\!\!\! Let $H=\left \langle a \right \rangle_2\ast\left \langle b^c \right \rangle_3$ be a subgroup of $G=G_1\ast G_2 = \left \langle a,b\ \!|\ \!a^2,b^3,[a,b] \right \rangle\ast\left \langle c|\ \!c^2 \right \rangle\simeq \mathbb{Z}_6\ast\mathbb{Z}_2$. Then $H$ is not a verbally closed subgroup in $G$ since for $f_1=ab\in G_1$ we have $(ab)^3 = a\in \left\langle a \right \rangle_2$ and $(ab)^2 = b^2\in \left \langle b \right \rangle_3$ but neither $(ab)^3 = 1$ nor $(ab)^2 = 1$.

{\Corollary \label{corom1} If $H$ is a verbally (algebraically) closed subgroup of $G=G_1\ast \cdots\ast G_n$, $1<\#G_i<\infty$, $i=1,\dots,n$, then it is a finitely generated group.}
\begin{Proof} Let $H$ be a verbally closed subgroup of $G$ with decomposition (\ref{Dec}). Then $F=\left \{ 1 \right \}$ by the first assertion of Theorem~\ref{ThM}. Since all free factors $H_{i,j}^{g_{i,j}}$ in this decomposition are finite groups we need to prove that there are only finitely many of them. If we assume the contrary, then (since the set of all subgroups of all group $G_i$, $i=1,\dots,n$ is finite) there are indices $i$, $j_1$ and $j_2$ such that $H_{i,j_1} = H_{i,j_2}$. But it means (by the second assertion of Theorem \ref{ThM}) that $H$ cannot be verbally closed in $G$.
\end{Proof}

Applying this corollary and the fact that a finitely generated algebraically closed subgroup of a finitely presented group is a retract (\cite{Romonkov}) we can conclude the following.

{\Corollary \label{corom2} For subgroups of $G=G_1\ast \cdots\ast G_n$, $1<\#G_i<\infty$, $i=1,\dots,n$ the properties of being algebraically closed and being a retract are equivalent.}

The conditions 1. and 2. of the theorem are necessary for verbal closedness but not sufficient. Indeed, it is clear that they hold for subgroup $H = \left \langle a \right \rangle_n\ast\left \langle b^c \right \rangle_n$ of $G = \left \langle a,b\ \!|\ \!a^n, b^n, [a,b] \right \rangle\ast\left \langle c\ \!|\ \!c^n \right \rangle\simeq \left ( \mathbb{Z}_n\times \mathbb{Z}_n \right )\ast \mathbb{Z}_n$, $n\ge 2$ or $n=\infty$. The goal of the two following theorems is to prove that if $n=2$ or $n=\infty$, then $H$ is not a verbally closed subgroup in $G$.

{\Theorem\label{ThEx} The subgroup $H = \left \langle a \right \rangle_2\ast\left \langle b^c \right \rangle_2$ of the group $G = \left \langle a,b\ \!|\ \!a^2, b^2, [a,b] \right \rangle\ast\left \langle c|c^2 \right \rangle\simeq \left ( \mathbb{Z}_2\times\mathbb{Z}_2 \right )\ast \mathbb{Z}_2$ is not verbally closed.
}

\begin{Proof} Since $\left \langle a,b\ \!|\ \!a^2, b^2, [a,b] \right \rangle\ast\left \langle c|c^2 \right \rangle\simeq \left \langle a,b,c\ \!|\ \!a^2,b^2,c^2,[a,b^c]\right \rangle$ it is sufficient to prove that $\widehat{H} = \left \langle a \right \rangle_2\ast\left \langle b \right \rangle_2$ is not a verbally closed subgroup in $\widehat{G} = \left \langle a,b,c\ \!|\ \!a^2,b^2,c^2,[a,b^c]\right \rangle.$

Let us consider the equation:

\begin{equation}\label{EqX}
\left(x^3\left[x,y^z\right]y^3\right)^2\left[x,y^z\right]^3=(ab)^2.
\end{equation}
This equation has a solution $x=a$, $y=b$, $z=c$ in $\widehat{G}$. We show that this equation has no solution in $\widehat{H} = \left \langle a \right \rangle_2\ast\left \langle b \right \rangle_2$.

Any element of $\widehat{H}$ can be expressed as $(ba)^k$ or $(ba)^ka$ where $k$ is an integer number. It means that, if we prove that any substitution of the form $x=(ba)^ka^{\varepsilon_1}$, $y=(ba)^ta^{\varepsilon_2}$, $z=(ba)^sa^{\varepsilon_3}$ where $k,t,s\in\mathbb{Z}$ and  $\varepsilon_1,\varepsilon_2,\varepsilon_3\in\{0,1\}\subset \mathbb{Z}$ is not a solution of equation~(\ref{EqX}), then we prove that this equation has no solution in $H$.
Consider the following eight possible cases and evaluate the left-hand side of the equation:

1) If $x=(ba)^k$, $y=(ba)^t$, $z=(ba)^s$, then $\left(x^3\left[x,y^z\right]y^3\right)^2\left[x,y^z\right]^3 = (ba)^{6(k+t)}$.

2) If $x=(ba)^ka$, $y=(ba)^t$, $z=(ba)^s$, then $\left(x^3\left[x,y^z\right]y^3\right)^2\left[x,y^z\right]^3 = (ba)^{-6t}.$

3) If $x=(ba)^k$, $y=(ba)^ta$, $z=(ba)^s$, then $\left(x^3\left[x,y^z\right]y^3\right)^2\left[x,y^z\right]^3 = (ba)^{6k}.$

4) If $x=(ba)^k$, $y=(ba)^t$, $z=(ba)^sa$, then $\left(x^3\left[x,y^z\right]y^3\right)^2\left[x,y^z\right]^3 = (ba)^{6(k+t)}.$

5) If $x=(ba)^ka$, $y=(ba)^ta$, $z=(ba)^s$, then $\left(x^3\left[x,y^z\right]y^3\right)^2\left[x,y^z\right]^3 = (ba)^{4(-k+t-s)}.$

6) If $x=(ba)^ka$, $y=(ba)^t$, $z=(ba)^sa$, then $\left(x^3\left[x,y^z\right]y^3\right)^2\left[x,y^z\right]^3 = (ba)^{6t}.$

7) If $x=(ba)^k$, $y=(ba)^ta$, $z=(ba)^sa$, then $\left(x^3\left[x,y^z\right]y^3\right)^2\left[x,y^z\right]^3 = (ba)^{6k}.$

8) If $x=(ba)^ka$, $y=(ba)^ta$, $z=(ba)^sa$, then $\left(x^3\left[x,y^z\right]y^3\right)^2\left[x,y^z\right]^3 = (ba)^{4(-k+s)}.$

It is clear that for any $k$, $t$ and $s$ the value of the left-hand side cannot be equal to $(ab)^2$ in $\widehat{H}$.
\end{Proof}

{\Theorem\label{ThLee1} The subgroup $H = \left \langle a \right \rangle_\infty\ast\left \langle b^c \right \rangle_\infty$ of the group $G = \left \langle a,b\ \!|\ \![a,b] \right \rangle\ast\left \langle c\ \!|\ \!\varnothing \right \rangle\simeq \left ( \mathbb{Z}\times \mathbb{Z} \right )\ast\mathbb{Z}$ is not verbally closed.
}

\begin{Proof} A word $w_r(x_1,\dots,x_m)$ is called a \emph{C-test word} in $m$ letters for the free group $F_r$ of rank $r\geqslant 2$ if for any two tuples $(g_1,\dots,g_m)$ and $(g_1',\dots,g_m')$ of elements of $F_r$ the following holds: if $w_r(g_1,\dots,g_m) = w_r(g_1',\dots,g_m')\ne 1$ in $F_r$, then there is an element $s\in F_r$ such that $g_i' = g_i^s,$ $i = 1,\dots,m$. In~\cite{Ivanov} Ivanov introduced and constructed first $C$-test words for $F_r$ in $m$ letters for any $r,m\geqslant 2$. In~\cite{Lee} Lee constructed for each $r,m\geqslant 2$ a $C$-test word $w_r(x_1,\dots,x_m)$ for $F_r$ with the addition property that $w_r(g_1,\dots,g_m)=1$ if and only if the subgroup of $F_r$ generated by $g_1,\dots,g_m$ is cyclic. We refer to such words as \emph{Lee words}.

Since $G = \left \langle a,b\ \!|\ \![a,b] \right \rangle\ast\left \langle c\ \!|\ \!\varnothing \right \rangle\simeq \widehat{G} = \left \langle a,b,c\ \!|\ \![a,b^c]\right \rangle$ it is sufficient to prove that $\widehat{H} = \left \langle a \right \rangle_\infty\ast\left \langle b \right \rangle_\infty$ is not a verbally closed subgroup in $\widehat{G} = \left \langle a,b,c\ \!|\ \![a,b^c]\right \rangle$. For this purpose consider the equation:

\begin{equation}\label{Eq1}
w_2\left(x_1,x_2,\left[x_1,x_2^{x_3}\right]\right) = w_2(a,b,1),
\end{equation}
where $w_2(x_1,x_2,x_3)$ is a Lee word for $F_2$ in three letters. Equation~(\ref{Eq1}) has a solution $x_1=a$, $x_2=b$, $x_3=c$ in $\widehat{G}$. We show that this equation has no solution in $\widehat{H}= F_2(a,b)$.

First, we note that since subgroup $\left \langle a,b \right \rangle = F_2(a,b)$ is not cyclic by the property of Lee words we have $w_2(a,b,1) \ne 1$ in $\widehat{H}$. If $x_i=h_i$, $i=1,2,3$ is a solution for (\ref{Eq1}) in $\widehat{H}$, then (by the definition of Lee words) there is an element $s\in \widehat{H}$ such that $h_1 = a^s$, $h_2 = b^s$ and $\left[h_1,h_2^{h_3}\right] = 1^s$. It is clear that this three equalities cannot be satisfied simultaneously in $\widehat{H}$ for any $s.$
\end{Proof}

\section{The proof of the main result}

One can easily see that under homomorphisms images of verbally closed subgroups are not such in general. Indeed, consider the group $\left \langle a \right \rangle_2\ast \left \langle b \right \rangle_2\ast \left \langle c \right \rangle_2$ and its subgroup $H = \left \langle c \right \rangle_2$. Since $H$ is a free factor of $G$ it is verbally closed in $G$ (actually, it is a retract in $G$). Let $\overline{H}$ be the image of $H$ under the natural homomorphism onto a factor group $G/R$, where $R$ is the normal closure of an element $[a,b]c$. The equation $[x,y]\overline{c}=\overline{1}$ has a solution $x=\overline{a}$, $y=\overline{b}$ in $G/R$ ($\overline{a}$, $\overline{b}$ and $\overline{c}$ are the images of $a$, $b$ and $c$, respectively). But this equation has no solution in $\overline{H}$ since it is a cyclic group and $\overline{c}\ne \overline{1}$ in $\overline{H}$. But under some additional restrictions the verbal closedness of subgroups is preserved for their homomorphic images.

\begin{Lemma}[\cite{Romonkov2}]\label{Lem1} Let the kernel $K=\ker\varphi$ of an epimorphism $\varphi\colon G \to \overline{G}$ be a verbal subgroup in $G$. If $H$ is a verbally closed subgroup in $G$, then $\overline{H} = \varphi(H)$ is a verbally closed subgroup in $\overline{G}$.
\end{Lemma}

\begin{Proof} Suppose that an equation
\begin{equation}\label{E1}
w(\underline{x})\overline{h}^{-1}=\overline{1},
\end{equation}

where $\overline{h} \in \overline{H}$ has a solution $\underline{x}=\overline{\underline{g}}$ in $\overline{G}$ (recall, that $\underline{x}$ and $\overline{\underline{g}}$ are the abridged notations for tuples of the form $(x_1,\dots,x_s)$ and $(\overline{g}_1,\dots,\overline{g}_s)$, respectively). We show that equation~(\ref{E1}) has a solution in $\overline{H}$.

Let $\underline{g}$ be a preimage of the tuple $\overline{\underline{g}}$ and $h$ be a preimage of $\overline{h}$ in $G$. Then in $G$ we have the equality: $$w(\underline{g})h^{-1}=k,$$

where $k \in K$.

Suppose that $K$ is generated by a set $V=\{f_\alpha(g_1,\dots,g_{n_{\alpha}})=f_{\alpha}\left ( \underline{g_{\alpha}} \right )\ \!|\ \!g_1,\dots,g_{n_{\alpha}}\in G\}$. We fix an expression $k=k\left ( f_{\alpha_1}\left ( \underline{g_{\alpha_1}} \right ),\dots, f_{\alpha_m}\left ( \underline{g_{\alpha_m}} \right )\right )$, $f_{\alpha_i}\left ( \underline{g_{\alpha_i}} \right ) \in V$, $i=1,\dots,m$ of the element $k$ and set up the following equation:

$$k\left ( f_{\alpha_1}\left ( \underline{y_1} \right ),\dots, f_{\alpha_m}\left ( \underline{y_m} \right )\right )^{-1}w(\underline{x})h^{-1}=1.$$

This equation has a solution $\underline{x}=\underline{g}$, $\underline{y_i}=\underline{g_{\alpha_i}}$, $i=1,\dots,m$ in $G$. Since $H$ is a verbally closed subgroup in $G$ this equation has a solution $\underline{x}=\underline{h}$, $\underline{y_i}=\underline{h_{\alpha_i}}$, $i=1,\dots,m$ in $H$. Since $K$ is a verbal subgroup in $G$ we have $k\left ( f_{\alpha_1}\left ( \underline{h_{\alpha_1}} \right ),\dots, f_{\alpha_m}\left ( \underline{h_{\alpha_m}} \right )\right )\in K$. From this we can conclude that the tuple $\underline{\varphi (h)}$ is a solution of equation~(\ref{E1}) in $\overline{H}.$
\end{Proof}

{\Corollary \label{corol3} If $H$ is a verbally closed subgroup in $G$ and $\varphi$ is an automorphism of $G$, then $\varphi(H)$ is a verbally closed subgroup in $G$.}

\begin{Lemma}\label{aux1}If $H$ is a verbally closed subgroup in $G$ and $\widehat{H}$ is a verbally closed subgroup in $H$, then $\widehat{H}$ is a verbally closed subgroup in $G$, i.e. verbal closedness is transitive.\end{Lemma}

\begin{Proof} Let $w(\underline{x})=\widehat{h}$, $\widehat{h}\in\widehat{H}$  be an equation which has a solution in $G$. Since $\widehat{h}\in H$ and $H$ is a verbally closed subgroup in $G$ this equation has a solution in $H$. Since $\widehat{H}$ is a verbally closed subgroup in $H$ we can conclude that the equation has a solution in $\widehat{H}$.
\end{Proof}

\begin{Lemma}\label{aux4} If $H=H_1\ast H_2$ is a verbally closed subgroup in $G$, then $H_1$ and $H_2$ are verbally closed subgroups in $G$.\end{Lemma}

\begin{Proof} Subgroup $H_1$ is a retract of $H$ (since it is a free factor of $H$), thus it is a verbally closed subgroup in $H$. By Lemma~\ref{aux1} $H_1$ is a verbally closed subgroup in $G$ (the same arguments hold for $H_2$).
\end{Proof}

Notice that if subgroups $H_1$ and $H_2$ are verbally closed in $G$, then subgroup $H_1\ast H_2$ may not be verbally closed. Indeed, subgroups $H_1=\left \langle a \right \rangle_2$ and $H_2=\left \langle b \right \rangle_2$ are retracts in $\left \langle a,b,c\ \!|\ \!a^2,b^2,c^2,[a,b^c]\right \rangle$, but $\left \langle a \right \rangle_2\ast\left \langle b \right \rangle_2$ is not a verbally closed subgroup in $\left \langle a,b,c\ \!|\ \!a^2,b^2,c^2,[a,b^c]\right \rangle$ (see Theorem~\ref{ThEx} for a proof).

To prove the first assertion of Theorem~\ref{ThM} we need to establish the following lemma.

\begin{Lemma}\label{LemF}

If $F_r$ is the free subgroup of rank $r\geqslant 1$ in group $G=G_1\ast \cdots\ast G_n$, $1<\#G_i<\infty$, $i=1,\dots,n$, then $F_r$ is not verbally closed in $G$.

\end{Lemma}

\begin{Proof} For each group $G_i$, $i=1,\dots,n$ we fix a set $S_i=\left\{s^i_{1},\dots,s^i_{t_i}\right\}$ of its generators. In this case $S=\bigcup_{i=1}^n S_i = \left\{s_{1},\dots,s_{t}\right\}$, where $t=t_1+\cdots+t_n$, is a generating set of $G$.

If $r=1$, then $F_1=\left \langle f \right \rangle_\infty$ for some element $f\in G$. Fix an expression $f=s_{i_1}s_{i_2}\cdots s_{i_m}$ of the element $f$ where all $s_{i_j}$ are elements of $S.$ Let $p$ be a prime number such that $p>\max\left ( \left \{ \#G_1,\dots,\#G_n \right \} \right )$, then the equation:
$$x_1^px_2^p\cdots x_m^p=f$$

has a solution $x_j=s_{i_j}^{k_j}$, $j=1,\dots,m$ in $G$ (if $s_{i_j} \in G_l$, then $k_j$ is a natural number such that $k_jp\equiv 1 \pmod{\# G_l}$). Since the left-hand side of the equation has a form $f^{pq}$, $q\in\mathbb{Z}$ for any substitution $x_j = f^{n_j}$, $j=1,\dots,m$ we can conclude that this equation has no solution in $F_1=\left \langle f \right \rangle_\infty$.

If $r\geqslant 2$ (possibly $r=\infty$), then $F_r = \left \langle f \right \rangle_\infty\ast F_{r-1}$ for some element $f\in G$. If $F_r$ is a verbally closed subgroup in $G$, then (by Lemma~$\ref{aux4}$) $\left \langle f \right \rangle_\infty$ is a verbally closed subgroup in $G$. But this contradicts the result of the preceding paragraph.
\end{Proof}

To prove the second assertion of the Theorem~\ref{ThM} we need to establish the following lemma.

\begin{Lemma}\label{Lem22} Let $G$ be a group of the type $G_1\ast \cdots\ast G_n$, $1<\#G_i<\infty$, $i=1,\dots,n$. If $H_1$ and $H_2$ are subgroups of a free factor $G_i$ and there are natural numbers $k_1$, $k_2$ and element $f\in G_i$ such that $f^{k_1}\in H_1$, $f^{k_1}\ne 1$ and $f^{k_2}\in H_2$, $f^{k_2}\ne 1$, then a subgroup of the form $H_1\ast H_2^g$, $g\in G$ is not verbally closed in $G$.
\end{Lemma}

To prove this lemma we need some auxiliary facts.

The following lemma is well-known.

{\Lemma\label{TreeLemma} Each automorphism of a directed tree has either a fixed vertex or an invariant line. If an automorphism has no fixed edges, then it has only one invariant line.}

\begin{Proof}
Let $f$ be an automorphism of a directed tree $T$ which fixes no vertices. Let us assume
$$m=\inf_{x\in V(T)}\!\!\left ( \rho \left ( x,f\left ( x \right ) \right ) \right ),$$
where $\rho(x,y)$ is the standard distance between vertices $x$ and $y$ and $V(T)$ is the set of all vertices of the tree $T$.

Since the distance $\rho(x,f(x))$ is a natural number for all $x\in V(T)$ there is a vertex $y\in V(T)$ such that $\rho(y,f(y))=m$. Let $z$ be the middle point of the segment $[y,f(y)]$, then $f(z)$ is the middle point of the segment $[f(y),f^2(y)]$ and $$\rho(z,f(z))\leqslant\rho(z,f(y))+\rho(f(y),f(z))=m.$$

Due to our choice of the vertex $y$ this inequality must be an equality. But the equality $\rho(z,f(z))=m$ means that the point $f(y)$ is located between the points $z$ and $f(z)$, and therefore between the points $y$ and $f^2(y)$. It follows that all points $f^n(y)$, $n\in\mathbb{Z}$ are located at the same line (which is an invariant line for the automorphism $f$).

Now we give a sketch of the proof of the second part of this lemma. Suppose that an automorphism $f$ has two invariant lines $a$ and $b$ (we still surmise that $f$ has no fixed vertices). In theory, there are four possibilities: the intersection of this lines is the empty set, a vertex, a segment or a beam.

If the intersection of the lines $a$ and $b$ is the empty set, then there is the unique segment $c$ which connects this lines. When $f$ acts on $T$ the ends of this segment must remain on their invariant lines, but it means that each edge of the segment $c$ is fixed. This contradicts our assumption.

The same arguments hold {\it mutatis mutandis} for the remaining three cases.
\end{Proof}

{\Definition\ }\!\!\! Any element $A$ of group $G=G_1\ast \cdots\ast G_n$, $1<\#G_i<\infty$, $i=1,\dots,n$ can be expressed as the reduced product $A = g_1\cdots g_s$ where each element $g_i$ lies in some free factor $G_k$. The \emph{norm} (or \emph{syllable length}) $|A|$ of an element $A$ is the number of syllables in its reduced product.

{\Lemma\label{BFL} Let $A$ be a cyclically reduced word of infinite order in group $G_1\ast \cdots\ast G_n$, $1<\#G_i<\infty$, $i=1,\dots,n$ and $g\in G$ be such element that the elements $A$ and $A^g$ do not commute in $G$. If $D$ is a cyclically reduced word in the conjugacy class of an element $A^{N_1}A^{gN_2}$, where $N_1,N_2\in\mathbb{N}$, $N_1,N_2\geqslant 2$, then the following estimate holds:
\begin{align*}
|D| > (N_1+N_2-4)|A|
\end{align*}}

\begin{Proof}

It is well known (see, for instance,~\cite{JohnMeier}) that every free product $G=G_1\ast \cdots\ast G_n$ of groups $G_i$, $i=1,\dots,n$ can be realized as a group of symmetries on some tree $T$. Furthermore, the stabilizers of the vertices of this tree under the action of group $G_1\ast \cdots\ast G_n$ are conjugates of $G_i$ for some $i$ and stabilizers of the edges are trivial. We fix such representation of $G$.

Since $A$ is an element of infinite order, the automorphism of the tree $T$ which corresponds to the action of $A$ has no fixed points, therefore (by Lemma~\ref{TreeLemma}) this automorphism has the unique invariant line (the same is true for $A^g$). Hereinafter we identify the elements of $G$ and the automorphisms of the tree $T$ that correspond to the actions of this elements.

Let $\alpha$ be the invariant line of the element $A$ and $\beta$ be the invariant line of the element $A^g$. Notice, that under the action of the elements $A$ and $A^g$ on the tree $T$ the vertices of their invariant lines $\alpha$ and $\beta$ are moved at the same distance (which is $|A|$). We show that the length $|I|$ of the intersection $I = \alpha\cap\beta$ is less than $2|A|$. Indeed, if $|I|\geqslant2|A|$, then there is a vertex $u\in V(I)$ such that either $A^2A^{2g}\in \textrm{St}(u)$ ($\textrm{St}(w)$ is the stabilizer of a vertex $w$ under the action) or $A^{-2}A^{2g}\in \textrm{St}(u)$ (this is true due to the fact that elements $A$ and $A^g$ move their invariant lines at the same distance and a point $u$ can be selected so that $A^{2g}u \in V(I)$).
Consider, for example, the first case (the same arguments hold {\it mutatis mutandis} for the second case). If $A^{2}A^{2g}\in \textrm{St}(u)$, then $AA^{g}\in \textrm{St}(u)$ and $AA^{g}\in \textrm{St}(A^{g}u)$ (to see it easily look at the vertices $u$, $A^gu$ and $A^{2g}u$ in Figure~\ref{fig1}), but the last means that under the action of the element $AA^{g}$ the segment $[u,A^{g}u]$ is fixed. Since under the action the stabilizers of the edges are trivial the element $AA^{g}$ must be the identity of $G$ but it is the contrary to the condition of the lemma.

\begin{otherlanguage}{english}

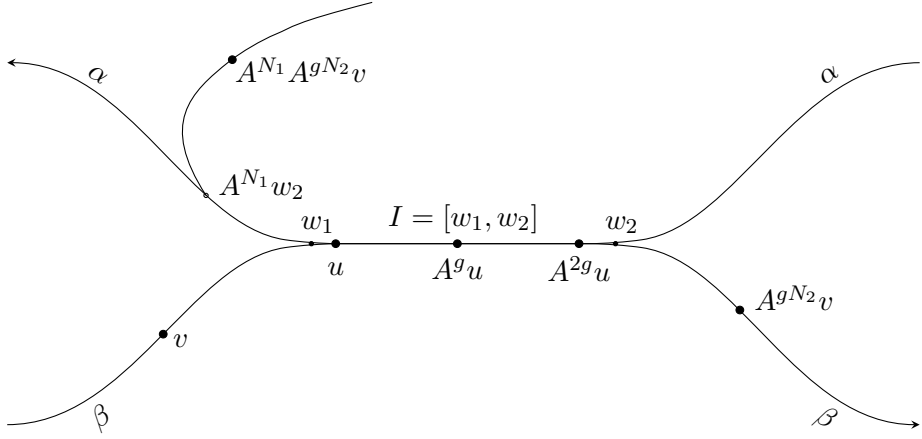
\begin{figure}[h]
\begin{center}
\begin{tikzpicture}[scale=0.8,>=stealth]

\filldraw [black] (-4.3,3.05) circle (1.8pt);
\filldraw [black] (-3,0) circle (1.0pt) node[xshift=0.5ex,yshift=1.5ex] {$w_1$};
\filldraw [black] (2,0) circle (1.0pt) node[xshift=0.5ex,yshift=1.5ex] {$w_2$};
\filldraw [black] (-2.6,0) circle (1.8pt) node[xshift=0.0ex,yshift=-2ex] {$u$};
\filldraw [black] (-0.6,0) circle (1.8pt) node[xshift=0.0ex,yshift=-2ex] {$A^g u$};
\filldraw [black] (1.4,0) circle (1.8pt) node[xshift=0.0ex,yshift=-2ex] {$A^{2g} u$};

\draw[rounded corners=10pt,name path=A,<-] (-8,3) .. controls (-6,3) and (-5,0.3) .. (-3,0) node[near start,sloped,above] {$\alpha$} -- (2,0) node[midway,xshift=0ex,yshift=2ex] {$I = [w_1,w_2]$}.. controls (4,0.3) and (5,3).. (7,3)

    node[near end,sloped,above] {$\alpha$};
\path [name path=T2] (-8,0.8) -- (0,0.8);
\draw [name intersections={of=A and T2, by=z}] (z) circle (1.pt) node[xshift=4.4ex,yshift=0.9ex] {$A^{N_1}w_2$};
\draw[rounded corners=10pt] (z) .. controls (-5.5,2) and (-5,2.5) .. (-4,3.3) node[xshift=4ex,yshift=-2ex] {$A^{N_1}A^{gN_2}v$}.. controls (-3.28,3.6) and (-3.35,3.65) ..(-2,4);

\draw[rounded corners=10pt,name path=B,->]  (-8,-3) .. controls (-6,-3) and (-5,-0.3) .. (-3,0) node[near start,sloped,below] {$\beta$} -- (2,0) .. controls (4,-0.3) and (5,-3).. (7,-3)

    node[near end,sloped,below] {$\beta$};

\path [name path=T] (-8,-1.5) -- (0,-1.5);
\filldraw [name intersections={of=B and T, by=x}] (x) circle (1.8pt) node[xshift=1.4ex,yshift=-0.5ex] {$v$};

\path [name path=T1] (8,-1.1) -- (0,-1.1);
\filldraw [name intersections={of=B and T1, by=y}] (y) circle (1.8pt) node[xshift=4.3ex,yshift=1ex] {$A^{gN_2}v$};

\end{tikzpicture}
\caption{Intersection of the invariant lines $\alpha$ and $\beta$ of the elements $A^{N_1}$ and $A^{gN_2}$. Arrows at the ends of lines indicate the directions in which the elements $A^{N_1}$ and $A^{gN_2}$ shift the lines.}
\label{fig1}
\end{center}
\end{figure}

\end{otherlanguage}

Now we can establish the estimate. To do it we establish the estimate on what minimum distance the element $A^{N_1}A^{gN_2}$ can move a vertex $v$ of the tree $T$. It is not difficult to see that the minimal shift is possible in the following case (this case is depicted in Figure~\ref{fig1}). Invariant lines $\alpha$ and $\beta$ intersect in a segment $I=[w_1,w_2]$ of the length $2|A|-1$ and vertices $v$ and $A^{gN_2}v$ are located on invariant line $\beta$ outside the segment $I$ and on the different sides with respect to this segment; the ends (vertices $w_2$ and $A^{N_1}w_2$) of the segments which connect the vertices  $A^{gN_2}v$ and $A^{N_1}A^{gN_2}v$ with the invariant line $\alpha$ and belong to this line lie on $\alpha$ on the different sides with respect to the segment $I$. In this case the distance between vertices $v$ and $A^{N_1}A^{gN_2}v$ is
 \begin{equation*}
\begin{split}
|[v,A^{N_1}A^{gN_2}v]| &= |[v,w_1]|+|[w_1,A^{N_1}w_2]|+|[A^{N_1}w_2,A^{N_1}A^{gN_2}v]|  =  \\
 & = |[v,A^{gN_2}v]|-|[w_2,A^{gN_2}v]| - |I| +|[w_2,A^{N_1}w_2]|-|I|+|[A^{N_1}w_2,A^{N_1}A^{gN_2}v]|=\\
 & = |[v,A^{gN_2}v]|+|[w_2,A^{N_1}w_2]|-2|I|\geqslant\\
 & \geqslant N_1|A|+N_2|A| - 2(2|A|-1) = (N_1+N_2-4)|A|+2,
\end{split}
\end{equation*}

where we used the following facts: $|[A^{N_1}w_2,A^{N_1}A^{gN_2}v]| = |[w_2,A^{gN_2}v]|$, $|[v,A^{gN_2}v]| = N_2|A|$ and $|[w_2,A^{N_1}w_2]| = N_1|A|$.

\end{Proof}

Now we can prove Lemma~\ref{Lem22}.

\begin{Proof} Let $H_1\ast H_2^g$ be a subgroup of $G$. For the elements $f$ and $g$ we fix expressions $f=f(s_1,\dots,s_{t}) = f(\underline{s})$ and $g=g(s_1,\dots,s_{t}) = g(\underline{s})$ where $s_i\in S$, $i=1,\dots,t$ ($S$ is the generation set of $G$ which has been described above during the proof of Lemma~\ref{LemF}). Let $N = 1+\prod_{i=1}^n\#G_i$, then the equation

\begin{equation}\label{E3}
f(x_1,\dots,x_t)^{k_1N}g(x_1,\dots,x_t)f(x_1,\dots,x_t)^{k_2N}g(x_1,\dots,x_t)^{-1} = f^{k_1}gf^{k_2}g^{-1}
\end{equation}
with variables $x_1,\dots,x_t$ has a solution $x_i=s_i$, $i=1,\dots,t$ in $G$. We denote $H_2^g$ as $H_2'$ and (with using abridged notations $f(x_1,\dots,x_t)$ = $f(\underline{x})$ and $g(x_1,\dots,x_t)$ = $g(\underline{x})$) rewrite equation~(\ref{E3}) in the following form:
\begin{equation}\label{E4}
f(\underline{x})^{k_1N}g(\underline{x})f(\underline{x})^{k_2N}g(\underline{x})^{-1} = ab,
\end{equation}
where $a=f^{k_1}\in H_1$, $a\ne 1$ and $b=gf^{k_2}g^{-1}\in H_2'$, $b\ne 1$. We show that equation~(\ref{E4}) has no solution in $H_1\ast H_2'$. Let a tuple $\underline{h}=(h_1,\dots,h_t)$ of elements of $H_1\ast H_2'$ be a solution of the equation. Thus, in $H_1\ast H_2'$ we have the following equality:
\begin{equation}\label{E5}
f(\underline{h})^{k_1N}g(\underline{h})f(\underline{h})^{k_2N}g(\underline{h})^{-1} = ab.
\end{equation}

We consider three cases.

1) Suppose that the element $f(\underline{h})$ has a finite order in $H_1\ast H_2'$. Then either $f(\underline{h})\in H_1^s$ or $f(\underline{h})\in H_2'^s$ for some element $s\in H_1\ast H_2'$. To be more specific, let $f(\underline{h})$ be an element of $H_1^s$, then the image of the equality~(\ref{E5}) under the natural epimorphism $\varphi\colon H_1\ast H_2'\to H_1\ast H_2'/[H_1,H_2']\simeq H_1\times H_2'$ has the following form:
$$\overline{f}^{k_1}\overline{h}_1\overline{f}^{k_2}\overline{h}_1^{-1} = \overline{a}\overline{b},$$
where $\overline{f},\overline{h}_1\in H_1$. But the last equality cannot hold in $H_1\times  H_2'$ since $\overline{f}^{k_1}\overline{h}_1\overline{f}^{k_2}\overline{h}_1^{-1}\in H_1$ and $\overline{a}\overline{b}\notin H_1$.

2) Suppose that the elements $f(\underline{h})$ and $g(\underline{h})f(\underline{h})g(\underline{h})^{-1}$ have infinite order and commute in $H_1\ast H_2'$. In this case there is an element $w\in H_1\ast H_2'$ such that $f(\underline{h}) = w^{m_1}$ and $g(\underline{h})f(\underline{h})g(\underline{h})^{-1} = w^{m_2}$ for some integer numbers $m_1$ and $m_2$ (indeed, in this case subgroup $\left \langle f(\underline{h}),g(\underline{h})f(\underline{h})g(\underline{h})^{-1} \right \rangle$ is an abelian group of infinite order and in the view of the Kurosh subgroup theorem it ought to have the form $\left \langle w \right \rangle_\infty$ for some element $w$). Thus, we can rewrite equality~(\ref{E5}) as:

$$w^{(k_1m_1+k_2m_2)N}=ab.$$

But this is impossible in $H_1\ast  H_2'$ since $N\geqslant 2$ and the element $ab$ is not a proper power in $H_1\ast H_2'$.

3) Suppose that the elements $f(\underline{h})$ and $g(\underline{h})f(\underline{h})g(\underline{h})^{-1}$ have infinite order and do not commute in $H_1\ast H_2'$. In this case there are natural numbers $m_1$, $m_2$ and a cyclically reduced word $A\in H_1\ast H_2'$ such that $f(\underline{h}) = A^{f_1m_1}$ and $g(\underline{h})f(\underline{h})g(\underline{h})^{-1}=A^{f_2m_2}$ for some elements $f_1,f_2\in H_1\ast H_2'$. Then it is straightforward to see that in this case the equality~(\ref{E5}) can be rewritten as:

\begin{align}\label{E7}
A^{k_1m_1N}A^{f_1^{-1}f_2k_2m_2N}=(ab)^{f_1^{-1}},
\end{align}

where the elements $A^{k_1m_1N}$ and $A^{f_1^{-1}f_2k_2m_2N}$ (and therefore the elements $A$ and $A^{f_1^{-1}f_2}$) do not commute in $H_1\ast H_2'$.

By Lemma~\ref{BFL} the left-hand side of the equality (\ref{E7}) is an element in the conjugacy class of a cyclically reduced word with the norm greater than $(k_1m_1N+k_2m_2N-4)|A|$. Since the right-hand side of this equality is the cyclically reduced word with the norm two and $((k_1m_1+k_1m_2)N-4)|A|\geqslant ((1\cdot1+1\cdot1)(1+4)-4)\cdot2=12$ (in this evaluation we used the facts: $\# G_i>1$ for all $i$, $n>1$, $k_1,k_2,m_1,m_2\in\mathbb{N}$ and $N = 1+\prod_{i=1}^n\#G_i$) we obtained the contradiction.
\end{Proof}

Now we are ready to prove Theorem~\ref{ThM}.

\begin{Proof} Let $H$ be a verbally closed subgroup in $G=G_1\ast \cdots\ast G_n$, $1<\#G_i<\infty$, $i=1,\dots,n$ with the Kurosh decomposition (\ref{Dec}).

1) If $H$ is a verbally closed subgroup in $G$, then (by Lemma~\ref{aux4}) each of its free factors is a verbally closed subgroup in $G$. In view of Lemma~\ref{LemF} a free subgroup can be verbally closed in $G$ if and only if it is the trivial group.

2) Suppose that for some index $i$ in decomposition~(\ref{Dec}) there are subgroups $H_{i,j_1}^{g_{i,j_1}}$ and $H_{i,j_2}^{g_{i,j_2}}$ and an element $f_i\in G_i$ such that $f_i^{k_1}\in H_{i,j_1}$, $f_i^{k_1}\ne 1$ and $f_i^{k_2}\in H_{i,j_2}^{g_i}$, $f_i^{k_2}\ne 1$ for some element $g_i\in G_i$ and natural numbers $k_1$ and $k_2$. By Lemma~\ref{aux4} subgroup $H_{i,j_1}^{g_{i,j_1}}\ast H_{i,j_2}^{g_{i,j_2}}$ is verbally closed in $G$, therefore (in accordance with Corollary~\ref{corol3}) subgroup $H_{i,j_1}\ast H_{i,j_2}^{g_{i,j_1}^{-1}g_{i,j_2}}$ is also verbally closed in $G$. But this is the contradiction to Lemma~\ref{Lem22}.
\end{Proof}

\newpage

\section{Acknowledgments}

The work of the author was supported by the Russian Foundation for Basic Research,
project no.~15-01-05823.

The author thanks A.~A.~Klyachko for many useful conversations and many useful remarks. Also the author is  grateful to the anonymous referee for useful remarks.

\begin{otherlanguage}{english}

\end{otherlanguage}

\end{document}